\documentclass[letterpaper, 10pt,reqno]{amsart}

\usepackage[usenames,dvipsnames]{xcolor}

\usepackage[portrait,margin=3cm]{geometry}
\usepackage{mathrsfs}

\usepackage{amssymb,amsthm,amsfonts,amsbsy,amstext,latexsym}
\usepackage[reqno]{amsmath}

\usepackage{color}
\usepackage{graphicx}
\usepackage{bbm}
\usepackage{comment}
\usepackage{a4wide,graphicx, listings,lscape, epstopdf, units, textcomp, fancyhdr, url, soul, color}
\usepackage{tikz}
\usetikzlibrary{shapes.geometric}
\usepackage[textsize=small]{todonotes}
\usepackage{enumitem}

\usepackage{mathtools}
\mathtoolsset{showonlyrefs}

\definecolor{MyDarkblue}{rgb}{0,0.08,0.50}
\definecolor{Brickred}{rgb}{0.65,0.08,0}

\usepackage{hyperref,upref}
\hypersetup{
	colorlinks=true,       % false: boxed links; true: colored links
	linkcolor=MyDarkblue,          % color of internal links
	citecolor=Brickred,        % color of links to bibliography
	filecolor=red,      % color of file links
	urlcolor=cyan           % color of external links
}

\newtheorem*{theorem*}{Theorem}
\newtheorem{theorem}{Theorem}[section]
\newtheorem{lemma}[theorem]{Lemma}

\theoremstyle{definition}

\newtheorem{remark}[theorem]{Remark}

\renewcommand{\P}{\mathbb{P}}
\newcommand{\prob}{\mathbb{P}}
\newcommand{\Pv}{\mathbb{P}}
\newcommand{\invisible}[1]{}
\newcommand{\EE}{\mathcal{E}}

\newcommand{\eps}{\varepsilon}

\newcommand{\cC}{\mathcal{C}}
\newcommand{\cE}{\mathcal{E}}

\newcommand{\cK}{\mathcal{K}}
\newcommand{\cO}{\mathcal{O}}

\newcommand{\eqn}[1]{\begin{equation} #1 \end{equation}}
\newcommand{\eqan}[1]{\begin{align} #1 \end{align}}
\newcommand{\sss}{\scriptscriptstyle}

\newcommand{\e}{{\mathrm e}}

\newcommand{\nn}{\nonumber}

\DeclareSymbolFont{extraup}{U}{zavm}{m}{n}
\DeclareMathSymbol{\varheart}{\mathalpha}{extraup}{86}
\DeclareMathSymbol{\vardiamond}{\mathalpha}{extraup}{87}

\newcommand{\ensymboldremark}{\hfill$\blacktriangleleft$}
\newcommand{\R}{\mathbb{R}}
\newcommand{\N}{\mathbb{N}}
\newcommand{\Z}{\mathbb{Z}}

\renewcommand{\emptyset}{\varnothing}

\newcommand*{\wt}{\widetilde}

\newcommand*{\be}{\begin{equation}}
	\newcommand*{\ee}{\end{equation}}
\newcommand*{\ba}{\begin{aligned}}
	\newcommand*{\ea}{\end{aligned}}
\newcommand*{\barr}{\begin{array}{c}}
	\newcommand*{\earr}{\end{array}}
\def \toinp    {\buildrel {\Pv}\over{\longrightarrow}}
\def \toindis  {\buildrel {d}\over{\longrightarrow}}

\makeatletter
\def\namedlabel#1#2{\begingroup
	#2%
	\def\@currentlabel{#2}%
	\phantomsection\label{#1}\endgroup
}

\makeatother
\newcommand{\bes}{\begin{equation*}}
	\newcommand{\ees}{\end{equation*}}

\renewcommand{\P}[1]{\mathbb{P}\!\left(#1\right)}
\newcommand{\E}[1]{\mathbb{E}\left[#1\right]}

\renewcommand{\N}{\mathbb{N}}

\renewcommand{\d}{\mathrm{d}}

\numberwithin{equation}{section}

\renewcommand{\e}{\mathrm{e}}

\newcommand{\T}{\mathbb T}
\newcommand{\D}[1]{\D_{#1}^{\sss (u,v)}}

\makeatletter
\newcommand{\leqnomode}{\tagsleft@true\let\veqno\@@leqno}
\newcommand{\reqnomode}{\tagsleft@false\let\veqno\@@eqno}
\makeatother
\newlength{\tagmarginsep} % Distance required
\tikzstyle{vertex}=[circle,fill=orange!60,minimum size=10pt,inner sep=0pt]
\tikzstyle{tedge} = [draw,ultra thick,->,>=stealth, orange]
\tikzstyle{esq}=[circle,fill=white,minimum size=10pt,inner sep=0pt]
\tikzstyle{up}=[<-,>=stealth]

\begin{document}
	
	\title[Long-range FPP on the complete graph]{Long-range first-passage percolation\\ on the torus}

        \author[VanderHofstad]{Remco van der Hofstad}
	\address{Department of Mathematics and Computer Science, Eindhoven University of Technology}
	\email{r.w.v.d.hofstad@tue.nl}
        \author[Lodewijks]{Bas Lodewijks}
	\address{Department of Mathematics, University of Augsburg}
	\email{bas.lodewijks@uni-a.de}
        \date{November $28$, $2023$}
	\maketitle
	
	\begin{abstract}
	    We study a geometric version of first-passage percolation on the complete graph, known as long-range first-passage percolation. Here, the vertices of the complete graph $\cK_n$ are embedded in the $d$-dimensional torus $\T_n^d$, and each edge $e$ is assigned an independent transmission time $T_e=\|e\|_{\T_n^d}^\alpha E_e$, where $E_e$ is a rate-one exponential random variable associated with the edge $e$, $\|\cdot\|_{\T_n^d}$ denotes the torus-norm, and $\alpha\geq0$ is a parameter. We are interested in the case $\alpha\in[0,d)$, which corresponds to the instantaneous percolation regime for long-range first-passage percolation on $\Z^d$ studied by Chatterjee and Dey~\cite{ChaDey16}, and which extends first-passage percolation on the complete graph (the $\alpha=0$ case) studied by Janson~\cite{Jans99}. We consider the typical distance, flooding time, and diameter of the model. Our results show a $1,2,3$-type result, akin to first-passage percolation on the complete graph as shown by Janson. The results also provide a quantitative perspective to the qualitative results observed by Chatterjee and Dey on $\Z^d$.    
	\end{abstract}

	
	\section{Introduction and main results}
	\label{sec-intro-main}

    The study of first-passage percolation was initiated by Broadbent and Hammersley~\cite{BroaHam57} and has since attracted an enormous wealth of attention and interest, first on infinite graphs such as $\Z^d$, and later also on finite (random) graphs. In this paper, we study \emph{long-range first-passage percolation on the complete graph}; a geometric analogue of first-passage percolation on the complete graph.

    First-passage percolation on the complete graph was first studied by Janson~\cite{Jans99}, where the edge-weights are independent and indentically distributed exponential random variables. Janson identified the leading asymptotics for the edge-weighted distance between fixed vertices, the flooding time (the maximum weighted distance over all pairs of vertices $(u,v)$ with $u$ fixed), and the diameter (the maximum weighted distance over all pairs of vertices), as well as more precise second-order behaviour. Later, this work was extended in, among others, \cite{AddBrouLug10,BhaHof17,HofHoogMie06,HoogMie08}. For the study of non-exponential edges, we refer to~\cite{BhaHof12,BhaHofHoog13,EckGoodHofNar13,EckGoodHofNar15.1,EckGoodHofNar15.2}, where first-passage percolation on the complete graph with edge-weights of the form $E^{s_n}$ is analysed. Here, $E$ is an exponential random variable and $(s_n)_{n\in\N}$ is a sequence that (possibly) depends on the number of vertices $n$. We note that these articles generalise some of the work for exponential edge-weights, where $s_n=1$ for all $n\in\N$.

    There is also a large body of work on first-passage percolation on \emph{random graphs}. Results for the Erd{\H o}s-R\'enyi random graph consider, among others, the typical distance, hopcount (the number of edges on the shortest edge-weighted path between two vertices), and the flooding time~\cite{BhaHofHoog11,HofHoogMie01,HofHoogMieg02}.

    In many \emph{scale-free} graphs, there is a clear distinction between graphs with degree distributions with finite- and infinite-variance degrees in terms of the behaviour of first-passage percolation. When the degrees have infinite variance, many such graph models are \emph{explosive} in that the typical weighted distance between vertices converges in distribution to a finite random variable.  The term \emph{explosion} is taken from the closely related theory of age-dependent branching processes. When this condition is not met, or the degrees have finite variance, first-passage percolation is \emph{conservative}, in the sense that the weighted distances diverge with the number of vertices in the graph.
    
    In the configuration model, the explosive and conservative settings have both been studied in detail in, among others, \cite{AdrKom18,BarHofKom17,BarHofKom19,BhaHofHoog10.2,BhaHofHoog17}, as well as in preferential attachment models~\cite{JorKom20,JorKom22}. Finally, in spatial models such as the scale-free percolation model and geometric inhomogeneous random graph model, both the conservative as well as the explosive setting have been studied~\cite{HofKom17,KomLod20}, too.

    In this paper, we study first-passage percolation on the complete graph, where the vertices are embedded on a $d$-dimensional torus $\T_n^d$ of volume $n$, and where each edge $e$ has an edge-weight $T_e:=\|e\|_{\T_n^d}^\alpha E_e$, where $\alpha\geq 0$ is a model parameter and the $E_e$ is a rate-one exponential random variable. One could consider this as first-passage percolation with a penalisation of long edges. The model is known as \emph{long-range first-passage percolation}. An infinite-size version of this model, essentially a `complete graph on $\Z^d$', has been studied by Chatterjee and Dey in \cite{ChaDey16}, who identified multiple phase transitions in the long-range parameter $\alpha$ for the behaviour of first-passage percolation. Related to this model, a degree-penalised version of spatial scale-free random graphs has recently been studied by Komj\'athy and co-authors in~\cite{KomLapLen21,KomLapLenSchal23_1,KomLapLenSchal23_2}, where multiple phases are discussed based on the penalisation strength.

    \paragraph{\bf Main contributions.} We study the typical distance, flooding time, and diameter of long-range first-passage percolation when the spatial dependence is weak, i.e.\  when $\alpha<d$, which Chatterjee and Dey call the \emph{instantaneous percolation regime}. They prove the qualitative result that, from any fixed vertex $v\in\Z^d$, any other vertex in $\Z^d$ can be reached within distance $t$ almost surely, for any $t>0$~\cite{ChaDey16}. Here, we study the finite-size setting on the complete graph, where we provide a quantitative counterpart to their results for the three aforementioned quantities. In particular, our results extend those of Janson \cite{Jans99} on the complete graph, which is the particular case $\alpha=0$.

	\subsection{Model}
	\label{sec-model}
    In this section, we define the model and fix notation.

    \noindent
	\textbf{Notation.} For $n,d\in\N$, let $\T_n^d$ be the discrete $d$-dimensional torus of volume $n$. That is, $\T_n^d:=[-n^{1/d}/2, n^{1/d}/2]\cap \Z^d$. It will be convenient to consider $\T_n^d$ as an equivalence relation on $\Z^d$, where $x$ is equivalent to $y$ when $x-y\in n^{1/d}\Z^d$. Then, we define the torus  $\infty$-norm
	\be 
	\|u\|_{\T_n^d}:=\min \|x\|_{\Z^d},
	\ee
	where $x$ runs over all elements in $\Z^d$ that are equivalent to $u$ and $\|\cdot \|_{\Z^d}$ denotes the Euclidean infinity norm. In a similar way, we can then define the torus $p$-norm on $\T_n^d$ by
	\be 
	\|u\|_{\T_n^d,p}:=\min \|x\|_{\Z^d,p},
	\ee 
	where again $x$ runs over all elements in $\Z^d$ that are equivalent to $u$, and  
        $$\|x\|_{\Z^d,p}=\Big(\sum_{i=1}^d x_i^p\Big)^{1/p}$$ 
    is the $p$-norm on $\Z^d$. From now on, we omit the subscript $\T_n^d$ for the norm when the context is clear, since we always work with the torus $p$-norm. Furthermore, we do specify which $p$ norm we use, since the results and proofs holds for all such norms, possibly by varying constants that do not alter the results.
	
	Let $G_n$ denote the complete graph with vertex set $\T_n^d$, i.e.\  $G_n:=(\T_n^d,\EE_n)$ where $\EE_n:=\{\{u,v\}\colon u\neq v\in \T_n^d\}$, and we study long-range first-passage percolation on it: Fix $\alpha\geq 0$ and assign to each edge $e\in \EE_n$ an edge-weight $T_e:=\|e\|^\alpha E_e$, where $(E_e)_{e\in \EE_n}$ is a collection of i.i.d.\ rate-one exponential random variables.
 
\subsection{Main result}
\label{sec-main}
 We start by introducing the main objects of study. For two vertices $u,v$, let $\Pi_{u,v}$ denote the set of all self-avoiding paths from $u$ to $v$. We define the transmission time between $u$ and $v$ as 
	\be 
	X_{u,v}:=\min_{\pi\in \Pi_{u,v}}\sum_{e\in \pi}E_e.
	\ee 
	We also define 
	\be \label{eq:rn}
	R_n:=\sum_{\substack{u\in \T_n^d\\ u\neq \mathbf0}}\|u\|^{-\alpha}.
	\ee 
	
	Our main result is as follows:
	
	\begin{theorem}[1,2,3 times $\log{n}/R_n$ for long-range first-passage percolation on the torus]
		\label{thrm:fpp}
		Fix $d\in\N, \alpha\in[0,d)$. Let $U$ and $V$ be two distinct vertices in $\T_n^d$, sampled uniformly at random, and let $u\in \T_n^d$ be some fixed vertex. Then, 
		\eqan{
			X_{U,V}\frac{R_n}{\log n}&\overset{\mathbb P}{\longrightarrow} 1;\label{typical}\\
			\max_{v\in \T_n^d}X_{u,v}\frac{R_n}{\log n}&\overset{\mathbb P}{\longrightarrow} 2; \label{flooding}\\ 
			\max_{u,v\in\T_n^d}X_{u, v}\frac{R_n}{\log n}&\overset{\mathbb P}{\longrightarrow} 3;\label{diameter}
		}
		where~\eqref{flooding} holds for any $u\in \T_n^d$.
	\end{theorem} 
	
	\begin{remark}[Interpretation main result]
		\label{rem-interpretation-main}
		As our proof will indicate, Theorem \ref{thrm:fpp} can be interpreted as follows. The typical distances in \eqref{typical} originate from a {\em continuous-time branching process} approximation of the neighborhood growth. This approximation kicks in on scale $1/R_n$, i.e.\  after rescaling the edge-weights by a factor $1/R_n$. Further, the vertices realizing the maxima in \eqref{flooding} are {\em traps}, meaning that their closest neighbor is at distance $(\log{n})/R_n$ away. Thus, the $v$ that realizes the maximum in \eqref{flooding} is one of these traps, while the $u$ and $v$ realizing the maximum in \eqref{diameter} are {\em both} traps. Even though our proof does not follow these lines, this intuition is very useful. It would be of interest to try to turn the above intuition into a proof.
		\ensymboldremark
	\end{remark}
	
	\begin{remark}[First-passage percolation on the complete graph]
		\label{rem-FPP-complete}
		When $\alpha=0$, all edges have an exponential weight with rate one, so that the spatial embedding of the vertices does not play a role. In this setting, we recover first-passage percolation on the complete graph and our results are in line with those of Janson~\cite[Theorem $1.1$]{Jans99}. Our proof for the lower bound on the flooding time and diameter in \eqref{flooding} and \eqref{diameter} follow that of Janson \cite{Jans99}. All other bounds rely on a large-deviation analysis of the exploration processes from two sources.
		\ensymboldremark
	\end{remark} 

        \begin{remark}[More general distance functions] The results of Theorem~\ref{thrm:fpp} also hold when using a distance function $r\colon \R_+\to \R_+$ such that $r$ is (eventually) decreasing and      $$
            \int_1^{n^{1/d}} r(x)x^{d-1}\,\d x/\log n\to\infty
            $$ 
        as $n$ tends to infinity, where we assign an edge $e$ the weight $E_e/r(\|e\|)$. In particular, $r(x)=\ell(x)x^{-\alpha}$ with $\ell$ a slowly varying function (i.e.\  $\ell$ such that $\lim_{x\to\infty}\ell(cx)/\ell(x)=1$ for any $c>0$) and $\alpha\in[0,d)$.
        \ensymboldremark
        \end{remark}
	
	\subsection{Discussion and open problems}
	\label{sec-disc-op}
	In this section, we discuss the structure of the proof of the main result, discuss the quantity $R_n$, and state some open problems.
	
	\paragraph{\bf The proof.} Our proof is organised as follows:
    In Section \ref{sec-LD-upper-bounds}, we perform a large-deviation analysis for typical distances that allows us to prove the upper bounds in \eqref{typical}, \eqref{flooding} and \eqref{diameter}: In Section \ref{sec:birth}, we discuss the continuous-time branching process approximation of the neighborhood growth, Section \ref{sec-growth-rate-CTBP} discusses the growth rate of this exploration process, Section \ref{sec-LD-diameter} performs a large-deviation analysis for upper bound on the diameter in \eqref{diameter}, while Section \ref{sec-LD-upper-bounds-typical-flooding} extends this large-deviation analysis to prove upper bounds on the flooding time in \eqref{flooding} and on the typical distances in \eqref{typical}. Having established all upper bounds, we move to the lower bounds in Section \ref{sec-lb-distances}. We prove the lower bound on typical distances in \eqref{typical} in Section \ref{sec-typical}, the lower bound on the flooding time in \eqref{flooding} in Section~\ref{sec:floodlb}, and the lower bound on the diameter in \eqref{diameter} in Section \ref{sec-diameter-lb}.
	
	\paragraph{\bf The limiting constant.}
	The dependence of $R_n$, which appears in Theorem~\ref{thrm:fpp} and is defined in~\eqref{eq:rn}, on $n$ can be made explicit. With relative ease, one can show that there exists a limiting constant $R=R(d,p,\alpha):=\lim_{n\to\infty} R_n n^{-(1-\alpha/d)}$ when $\alpha<d$. This implies $R_n$ is of the order $n^{1-\alpha/d}$ and provides a quantitative limit for the results in Theorem~\ref{thrm:fpp}. An interesting question is what the value $R$ is for different choices of $d,p,$ and $\alpha$. 
	
	In Appendix~\ref{sec:app} we provide an integral representation for $R$ for any dimension and any norm. Moreover, we explicitly identify $R$ (or provide a recursive approach how one can do so) for $p=1$ and $p=\infty$ and all dimensions, as well as for dimension two and any $p\in[1,\infty)\cup\{\infty\}$.

    \paragraph{\bf Comparison to long-range percolation.}   In long-range percolation on $\Z^d$, there is a number of phase transitions for the typical distance (also referred to as chemical distance) between, say, the origin and a `far away' vertex $x$. In this model, non-nearest neighbour edges between vertices $x,y\in \Z^d$ are created independently with probability $1-\exp(\beta r(x-y))$, where $r(x)\sim |x|^{-s}$ as $|x|\to\infty$ for some $s>0$. Here, Benjamini et al.~\cite{BenKesPerSch11} show that the model has finite diameter $\lceil d/(d-s)\rceil$ when $s<d$. Our parameter $\alpha$ can be thought of as the parameter $s$ in the long-range parameter model. We choose to be consistent with the notation in the work of Chatterjee and Dey on long-range first-passage percolation in~\cite{ChaDey16} and thus use $\alpha$ rather than $s$. The results presented here are similar to the long-range percolation setting in the sense that the diameter of the torus is bounded, and in fact converges to zero with $n$. The latter is due to the fact that we consider an edge-weighted distance, where the random edge weights have support in $(0,\infty)$, so that optimal paths use edges with vanishing edge-weights, whereas the graph distances studied for long-range percolation are always at least one. We refer the interested reader to the article of Biskup and Krieger~\cite{BisKrie24} for a more in-depth overview of results for chemical distances in long-range percolation.
    
	\paragraph{\bf Open problems.} For the complete graph, i.e.\  for $\alpha=0$, Janson~\cite[Theorem $1.1$]{Jans99} discusses the {\em fluctuations} around the main asymptotics in \eqref{typical} and \eqref{flooding}, while the fluctuations around \eqref{diameter} were identified by Bhamidi and the first author in \cite{BhaHof17}. It would be of interest to extend that analysis to our spatial setting. Janson also considers the \emph{hopcount} $H_{u,v}$, i.e.\  the number of edges on the shortest edge-weighted path between $u$ and $v$, as well as $\max_{v\in \T_n^d} H_{u,v}$, in the complete graph. Addario-Berry, Broutin, and Lugosi consider $\max_{u,v\in\T_n^d} H_{u,v}$ in~\cite{AddBrouLug10}. It would be interesting to see if these results can be extended to the spatial setting.
 
    It would also be of interest to see whether the result on the typical distance, as in~\eqref{typical}, can be proved to hold for \emph{fixed} vertices $u$ and $v$, rather than vertices selected uniformly at random. Though we expect this to be the case, the proof of the \emph{lower bound} crucially depends on the uniformity of these two vertices (in fact, it suffices if one is uniformly chosen, by translation invariance). The \emph{upper bound} does apply in full generality. In the complete graph setting, i.e.\  for $\alpha=0$, independently of which vertices have already been explored, the next vertex that is to be explored is uniform among all unexplored vertices. This is only true in the spatial setting when the vertices $U$ and $V$ are uniform themselves. It seems to us that  understanding the distribution of the set of explored vertices, and how close it is to the uniform distribution (among all sets in $\T_n^d$ of size equal to the number of explored vertices), is important to enable one to prove~\eqref{typical} for fixed vertices. We expect this understanding to also be crucial in analysing the other statistics mentioned above.
 
    Further, it would be of interest to extend the results in this paper to other edge-weight distributions or transmission times. For example, Bhamidi and the first author in \cite{BhaHof12} consider the natural example of a {\em power} of an exponential distribution, and identify the continuous-time branching-process approximation of the first-passage percolation exploration, which might yield a starting point for the present spatial problem.

    The case $\alpha>d$ is studied in detail on the $d$-dimensional lattice by Chatterjee and Dey in~\cite{ChaDey16}. They determine three additional regimes (that is, $\alpha\in (d,2d)$, $\alpha\in(2d,2d+1)$, and $\alpha>2d+1$) in which the edge-weighted distance between the origin and a vertex $x\in\Z^d$ with $\|x\|=m$ grows as a function of $m$. Since the distance grows with $m$ in these regimes, we expect analogous behaviour for the typical distance on the torus that we study here, with $m$ replaced by $n^{1/d}$. The flooding time and diameter when $\alpha>d$ should be of the same order as the typical distance. Indeed, the time to escape traps, i.e.\ the vertices that realise the maxima in~\eqref{flooding} and~\eqref{diameter} (as mentioned in Remark~\ref{rem-interpretation-main}), are of order $\log n$. On the other hand, the typical distances grow much faster than $\log n$ when $\alpha>d$ (as proved by Chatterjee and Dey in~\cite{ChaDey16}).

	\section{Large deviations approach for distances: upper bounds}
	\label{sec-LD-upper-bounds}
	
	In this section we analyse our distances in $\T_n^d$ under the random metric. We start by proving the upper bound in \eqref{diameter} in  Theorem~\ref{thrm:fpp}. For this, we use a large deviations approach, similar to Janson~\cite{Jans99}. At the end of this section, we show how this argument can easily be extended to the upper bounds in \eqref{typical} and \eqref{flooding} in  Theorem~\ref{thrm:fpp}.
	
	\subsection{Long-range first-passage percolation as an exploration process}
	\label{sec:birth}
	
	We can view the exploration process of long-range first-passage percolation process on $G_n$, started from any vertex $u\in\T_n^d$, as a birth process $T_n(u,t)$. Any vertex in this birth process has a unique spatial location in $\T_n^d$. Here, $T_n(u,0)=1$ and for any $j\in[n-1]$ and $t\geq 0$, when $T_n(u,t)=j$, a birth takes place at rate 
	\be \label{eq:lambdaj}
	\Lambda_j:=\sum_{\ell=0}^{j-1}\bigg(R_n-\sum_{\substack{i=0\\ i\neq \ell}}^{j-1}\|v_i-v_\ell\|^{-\alpha}\bigg),
	\ee 
	where $v_i$ denotes the $i^{{\rm th}}$ vertex born in the birth process $T_n(u,t)$ and $v_0=u$. Here, the argument of the outer sum equals the rate at which the vertex $v_\ell$ gives birth to the next vertex. It follows that $\Lambda_j$ is random. When the $j^{{\rm th}}$ birth takes place, conditionally on $\Lambda_j$ and $v_1,\ldots, v_{j-1}$, 
    \be\label{eq:spreadprob}
    \P{v_j=z\,|\, \Lambda_j,v_1,\ldots,v_{j-1}}=\frac{\sum_{i=0}^{j-1}\|z-v_i\|^{-\alpha}}{\Lambda_j}, \qquad z\in \T_n^d\backslash\{v_0,\ldots, v_{j-1}\}. 
    \ee 
    Then, observe that, for each $\ell\in\{0,\ldots, j-1\}$, the vertices $v_0,\ldots v_{\ell-1},v_{\ell+1},\ldots, v_{j-1}$ are at least as far away from $v_\ell$ as the $j-1$ closest vertices to $v_\ell$. The sum of the weighted distances of these $j-1$ closest vertices equals $R_j$. Using this thus yields the bounds
    \be \label{eq:Lambdabounds}
    j(R_n-R_j)\leq \Lambda_j\leq j R_n.
    \ee 
    We thus find that for all $j=o(n)$ almost surely,  since $\alpha<d$,
        \eqn{
        \label{bounds-Gamma-j}
        jR_n(1-\cO((j/n)^{1-\alpha/d}))\leq \Lambda_j\leq jR_n.
        }
    Moreover, if we let $E_j\sim \text{Exp}(\Lambda_j)$, and, for $j\in[n-1]$,
	\eqn{
		\label{stochastic-bounds}
		E_{j,U}\sim \text{Exp}(j(R_n-\cO(j^{1-\alpha/d}))),
    \qquad E_{j,L}\sim \text{Exp}(jR_n),
	}
	where the $(E_{j,L})_{j\in[n-1]}$ (resp.\ $(E_{j,U})_{j\in[n-1]}$) are independent, then $E_{j,L}\preceq E_j\preceq E_{j,U}$ for any $j=o(n)$, where $X\preceq Y$ denotes that the random variable $X$ stochastically dominates the random variable $Y$. As a result, the time until $j=o(n)$ vertices have been discovered stochastically dominates, and is stochastically dominated by, a sequence of independent exponential random variables with deterministic rates that are also controllably close to one another.
	
	\subsection{Growth rate of exploration process}
    \label{sec-growth-rate-CTBP}
	
	We describe our first-passage percolation process started from a vertex $v\in\T_n^d$ via the random metric induced by the random edge-weights. We define 
	\be 
	B_v(t):=\{u\in \T_n^d\colon  X_{v,u}\leq t\}, 
	\ee 
	as the (random) ball of radius $t$ around $v$ with respect to the transmission-time metric. That is, $B_v(t)$ denotes all vertices in $T_n^d$ that can be reached within time $t$ from $v$. Let us also define $(\tau_k)_{k\in\N_0}$ by $\tau_k:=\inf\{t>\tau_{k-1}\colon  |B_v(t)|=k+1\}, k\in\N$, and $\tau_0=0$. It is clear that the intervals $\tau_k-\tau_{k-1}$ for $k\in\N$ are exponentially distributed with rate $\Lambda_k$, conditionally on $\Lambda_k$. Since $\Lambda_k\approx kR_n$ follows from~\eqref{bounds-Gamma-j} for $k$ not too large, we have the following lemma, which controls the fluctuations of the times at which the first-passage percolation process discovers vertices in $\T_n^d$:
    
    \begin{lemma}[Distributional limit of $\tau_k$]\label{lemma:tauk}
       Fix $\alpha\in(0,d)$. Let $k=k(n)$ be such that $k\to\infty$ and $k=o(n)$. Then,
        \be 
        R_n\tau_k- \log k\toindis \Lambda,
        \ee 
        where $\Lambda$ is a Gumbel random variable. In particular, when $k=n^\beta(1+o(1))$ for some $\beta\in(0,1)$, 
    \be 
    R_n\tau_k-\beta \log n\toindis \Lambda, 
    \qquad
    {\rm and~thus }
    \qquad 
    \label{eq:tauinp}
    \tau_k \frac{R_n}{\log n}\toinp \beta.
    \ee 
    \end{lemma}
    Lemma \ref{lemma:tauk} proves more than we will need below, since it also identifies the \emph{fluctuations} of $\tau_k$. One can interpret the Gumbel limiting variable $\Lambda$ as $\log(1/E)$, where $E$ is exponential with rate 1. This variable $E$ is the limit of ${\mathrm e}^{-t}Z_t$, where $(Z_t)_{t\geq 0}$ is the Yule process that arises as the distributional limit of the (rescaled) $(|B_v(t/R_n )|)_{t\geq 0}$ discovery process.
    
    \begin{proof} 
    Let $(E_k)_{k\in\N}$ be a sum of i.i.d.\ exponential random variables with rate one. It holds that $\tau_{k+1}$ is equal in distribution to
    \be 
    \tau_{k+1}\overset d= \sum_{j=1}^{k} \frac{E_j}{\Lambda_j}, 
    \ee 
    where we recall $\Lambda_j$ from~\eqref{eq:lambdaj}, the (random) rate at which the first-passage percolation process spreads from $j$ vertices to $j+1$ vertices. As a result, we can write 
    \be 
    R_n\tau_{k+1}-\log k=\sum_{j=1}^k \frac{E_j}{j}-\sum_{j=1}^k \frac1j +\gamma +\frac{1}{2k}-\eps_k+R_n\sum_{j=1}^k\Big(\frac{1}{\Lambda_j}-\frac{1}{jR_n}\Big)E_j, 
    \ee 
    where we use that $\log k=H_k+\gamma+1/(2k)-\eps_k$, where $H_k$ is the harmonic sum, $\gamma$ is the Euler-Mascheroni constant, and $\eps_k$ is an error term such that $\eps_k\in[0,1/(8k^2)]$. We now combine the first three terms, which contribute to the limiting random variable, and write all deterministic terms that tend to zero with $k$ as $o(1)$, to obtain 
    \be 
    \sum_{j=1}^k \frac{E_j-1}{j}+\gamma+\sum_{j=1}^k \frac{jR_n-\Lambda_j}{j\Lambda_j}E_j+o(1). 
    \ee 
    For the first sum, we have that as $k\to\infty$,
    \be 
    \sum_{j=1}^k \frac{E_j-1}{j}+\gamma\toindis \sum_{j=1}^\infty \frac{E_j-1}{j}+\gamma \overset d= \Lambda, 
    \ee 
    where $\Lambda$ is a Gumbel random variable. For the second sum, we use~\eqref{bounds-Gamma-j} to obtain 
    \be 
    jR_n-\Lambda_j\leq jR_{j-1}=\cO(j^{2-\alpha/d}).
    \ee 
    Moreover, using that $\Lambda_j=jR_n(1-o(1))$ almost surely for any $j=o(n)$, as follows from~\eqref{eq:Lambdabounds}, we obtain for some constant $C>0$, 
    \be
    \sum_{j=1}^k \frac{jR_n-\Lambda_j}{j\Lambda_j}E_j=\frac{C}{R_n}\sum_{j=1}^k j^{-\alpha/d}E_j.
    \ee 
    It thus remain to show that sum is smaller than $\eps R_n$ with high probability for any $\eps>0$ to conclude that the sum on the left-hand side converges to zero. So, using a Markov bound,
    \be 
    \P{\sum_{j=1}^k j^{-\alpha/d}E_j\geq \eps R_n}\leq \frac{1}{\eps R_n}\sum_{j=1}^k j^{-\alpha/d}=\cO(k^{1-\alpha/d}/R_n)=o(1), 
    \ee 
    since $k=o(n)$ (as $R_n=(R+o(1))n^{1-\alpha/d}$ when $\alpha<d$). So, for any $k$ that diverges with $n$ such that $k=o(n)$, 
    \be 
    R_n\tau_{k+1}-\log k\toindis \Lambda.
    \ee
    The desired result then follows since $\log k-\log(k+1) =o(1)$.
    \end{proof}	
	
	\subsection{Large deviation analysis for upper bound diameter}
	\label{sec-LD-diameter}
	We have, for any $\eps>0$, 	
	\be\ba
	\P{\max_{u,v\in \T_n^d}X_{u,v}\geq (3+2\eps)R_n^{-1}\log n}
	&\leq \sum_{u,v\in \T_n^d}\P{ X_{u,v}\geq (3+2\eps)R_n^{-1}\log n}\\
	&\leq n^2\max_{u,v\in \T_n^d}\P{ X_{u,v}\geq (3+2\eps)R_n^{-1}\log n}.
	\ea\ee 
	We now bound the transmission time between $u$ and $v$ from above by 
	\be \label{eq:timeub}
	\tau_{a_n}^{\sss (u)}+\tau_{a_n}^{\sss (v)}+\min_{x\in B_u(\tau_{a_n}^{\sss (u)}), y\in B_v(\tau_{a_n}^{\sss (v)})}T_{xy}.
	\ee 
    Here, $\tau_{a_n}^{(v)}$ denotes $\tau_{a_n}$ for a FPP process started from $v\in\T_n^d$. Note that both $\tau_{a_n}^{(v)}$ and $\tau_{a_n}^{(u)}$ satisfy~\eqref{eq:tauinp} in Lemma~\ref{lemma:tauk} (so also jointly), but they are independent if $B_u(\tau_{a_n}^{\sss (u)})\cap B_v(\tau_{a_n}^{\sss (v)})=\emptyset$ holds.
	Indeed, if $B_u(\tau_{a_n}^{\sss (u)})\cap B_v(\tau_{a_n}^{\sss (v)})=\emptyset$ holds, then the time for the exploration processes of $u$ and $v$ to both reach size $a_n$ plus the shortest edge between their respective explored vertices, is a possible way to connect $u$ and $v$, and thus an upper bound for the transmission time. If, on the other hand, $B_u(\tau_{a_n}^{\sss (u)})\cap B_v(\tau_{a_n}^{\sss (v)})\neq \emptyset$, then the third term in~\eqref{eq:timeub} equals zero almost surely, and the transmission time between $u$ and $v$ is no more than the time it took for both exploration processes to reach size $a_n$. As such, we obtain the upper bound
	\be \ba \label{eq:nsquarebound}
	n^2{}&\max_{u,v\in \T_n^d}\P{ \tau_{a_n}^{\sss (u)}+\tau_{a_n}^{\sss (v)}+\min_{x\in B_u(\tau_{a_n}^{\sss (u)}), y\in B_v(\tau_{a_n}^{\sss (v)})} T_{xy}\geq (3+2\eps)R_n^{-1}\log n}\\
	&\leq n^2\!\max_{u,v\in \T_n^d}\bigg[\P{ \tau_{a_n}^{\sss (u)}+ \tau_{a_n}^{\sss (v)}\geq\frac{(3+\eps)\log n}{R_n}}+\P{\min_{x\in B_u(\tau_{a_n}^{\sss (u)}), y\in B_v(\tau_{a_n}^{\sss (v)})}\!\!\!\! T_{xy}\geq \frac{\eps \log n}{R_n}}\bigg].
	\ea \ee 
	Let us bound the second probability on the right-hand side first. The minimum contains $(a_n+1)^2$ many terms. As such, by bounding the spatial distance between any $x\in B_u(\tau^{\sss (u)}_{a_n})$ and $y\in B_v(\tau_{a_n}^{\sss (v)})$ from above by $Cn^{1/d}$ for some sufficiently large constant $C>0$, the minimum is stochastically dominated by an exponential random variable with rate $(a_n+1)^2(Cn^{1/d})^{-\alpha}=(c_1^2C^{-\alpha}+o(1))n^{1-\alpha/d}$. Hence,   
	\be \label{eq:bridge}
	\P{\min_{x\in B_u(\tau_{a_n}^{\sss (u)}), y\in B_v(\tau_{a_n}^{\sss (v)})} T_{xy}\geq \eps R_n^{-1} \log n}\leq \exp(-(c_1^2C^{-\alpha}R^{-1}\eps+o(1)) \log n)=o(n^{-2}), 
	\ee  
	where the final equality holds when $c_1>\sqrt{2RC^\alpha/(\eps(1-\alpha/d))}$. Observe that this upper bound is independent of the choice of $u,v\in\T_n^d$. It thus remains to upper bound
	\be 
	\P{\tau_{a_n}^{\sss (u)}+ \tau_{a_n}^{\sss (v)}\geq(3+\eps)R_n^{-1}\log n}.
	\ee  
	Again, we aim to obtain an upper bound that is $o(n^{-2})$, independently of $u$ and $v$, so that combining this with~\eqref{eq:bridge} in~\eqref{eq:nsquarebound} yields that $\max_{u,v\in\T_n^d}X_{u,v}$ is at most $(3+2\eps)R_n^{-1}\log n$ with high probability, for any fixed $\eps>0$. To this end, we replace $\tau_{a_n}^{\sss (v)}$ by $\widehat\tau_{a_n}^{\sss (v)}$, which is the time it takes to reach $a_n$ vertices from $v$ when all edges incident to the vertices in $B_u(\tau_{a_n}^{\sss (u)})$ are removed. As this leaves fewer possible edges to reach vertices (and fewer vertices to be reached), $ \tau_{a_n}^{\sss (v)}\preceq \widehat\tau_{a_n}^{\sss (v)}$. Let us write $(z_j)_{j=0}^{a_n}$ for the vertices that are explored by $u$, in chronological order. That is, $z_i$ is the $i^{\text{th}}$ vertex explored by $u$ (with $z_0=u$): $X_{u,z_j}>X_{u,z_i}$ for $i\in[j-1]$ and $X_{u,z_j}<X_{u,z_k}$ for any $k>j$. Then, conditionally on $B_u(\tau_{a_n}^{\sss (u)})$, we know that $\widehat \tau_{a_n}^{\sss (v)}$ is a sum of $a_n$ exponential random variables $(\wt E_j)_{j\in[a_n]}$, where 
	\be 
	\wt E_j\sim \text{Exp}\bigg(\sum_{\ell=0}^{j-1}\bigg(R_n-\sum_{\substack{i=0\\ i\neq \ell}}^{j-1}\|z_i-z_\ell\|^{-\alpha}-\sum_{w\in B_u(\tau_{a_n}^{\sss (u)})}\|z_\ell-w\|^{-\alpha}\bigg)\bigg).
	\ee
	Similar to~\eqref{eq:Lambdabounds}, for each $j\in[a_n]$ we bound this rate from below by considering the $2a_n+1$ closest vertices to $z_\ell$ for each $\ell\in\{0,\ldots, j-1\}$. This yields a sequence of independent exponential random variables $(\widehat E_j)_{j\in[a_n]}$ such that $\wt E_j\preceq \widehat E_j$ for each $j\in[a_n]$, and 
	\be 
	\widehat E_j\sim \text{Exp}\big(j(R_n-R_{2a_n+1})\big), \qquad j\in[a_n-1].
	\ee 
	By~\eqref{bounds-Gamma-j} and since $2a_n+1=o(n)$, it follows that the rate equals $jR_n(1-o(1))$, where the $o(1)$ is independent of $j$. We can thus couple $\widehat \tau_{a_n}^{\sss (v)}$ with the $(\widehat E_j)_{j\in [a_n]}$ such that 
	\be 
	\widehat \tau_{a_n}^{\sss (v)}\preceq \sum_{j=1}^{a_n}\widehat E_j.
	\ee 
	Similarly, as discussed in Section~\ref{sec:birth}, we can stochastically bound $\tau_{a_n}^{\sss (u)}$ from above by the sum of the random variables $(E_{j,U})_{j\in[a_n]}$, which are independent of the $(\widehat E_j)_{j\in[a_n]}$. Combined with a Chernoff bound with $s=cR_n$ and $c\in(2/(2+\eps),1)$, this then yields
	\be \ba
	\P{\tau_{a_n}^{\sss (u)}+ \tau_{a_n}^{\sss (v)}\geq(3+\eps)R_n^{-1}\log n}&\leq n^{-(3+\eps)c}\prod_{j=1}^{a_n}\E{\e^{sE_{j,U}}}\prod_{j=1}^{a_n}\E{\e^{s\widehat E_j}}\\
	&=n^{-(3+\eps)c}\prod_{j=1}^{a_n}\Big(\frac{jR_n(1-o(1))}{jR_n(1-o(1))-s}\Big)^2\\
	&\leq n^{-c(3+\eps)}\exp\Big(2\sum_{j=1}^{a_n}\frac{c}{j(1-o(1))-c}\Big)\\
	&=\exp(-c(3+\eps)\log n +2c\log(a_n)+o(\log n))\\
	&=n^{-c(2+\eps)+o(1)}.
	\ea \ee 
	By the choice of $c$, this upper bound is $o(n^{-2})$. Combining this with~\eqref{eq:bridge} in~\eqref{eq:nsquarebound} finally yields the desired result. This proves the upper bound on \eqref{diameter} in Theorem~\ref{thrm:fpp}.
	\qed
	
	\subsection{Large deviation analysis for upper bounds flooding time and typical distances}
	\label{sec-LD-upper-bounds-typical-flooding}
	We next extend this analysis to upper bounds on the flooding time in \eqref{flooding}, and on the typical distances in \eqref{typical}.
	
	\paragraph{\bf Upper bound on flooding time.}
    A similar approach, now taking $c\in (1/(1+\eps),1)$ and $u\in \T_n^d$ fixed, gives 
	\be \ba 
	\P{\max_{\substack{v\in \T_n^d\\ v\neq u}} X_{u,v}\geq (2+2\eps)R_n^{-1}\log n}&\leq n\max_{\substack{v\in \T_n^d\\ v\neq u}}\P{ X_{u,v}\geq (2+2\eps)R_n^{-1}\log n}\\
	&\leq \mathcal O(n^{1-(1+\eps)c+o(1)})+o(n^{-1})=o(1).
	\ea\ee 
    This proves the upper bound on~\eqref{flooding} in Theorem~\ref{thrm:fpp}.
	\qed
	
	\paragraph{\bf Upper bound on typical distance.}
	Again using a similar approach, for distinct $u,v\in \T_n^d$ and with $c\in(0,1)$, 
	\be \ba
	\P{X_{u,v}\geq (1+2\eps)R_n^{-1}\log n}&\leq \P{\tau_{a_n}^{\sss (u)}+\tau_{a_n}^{\sss (v)}+\min_{x\in B_u(\tau_{a_n}^{\sss (u)}), y\in B_v(\tau_{a_n}^{\sss (v)})}T_{uv}\geq (1+2\eps)R_n^{-1}\log n}\\
	&= \cO(n^{-c\eps}).
	\ea \ee 
	As this holds for \emph{any} distinct $u,v\in\T_n^d$ and the upper bound is independent of the choice of $u$ and $v$, this also extends to the typical distance between $U$ and $V$, two distinct vertices selected uniformly at random. This proves the upper bound on \eqref{typical} in Theorem~\ref{thrm:fpp}. In particular, while we state \eqref{typical} in Theorem~\ref{thrm:fpp} only for two uniform vertices, the upper bound actually holds for {\em all} fixed pairs as well.
	\qed
	
	\section{Lower bounds on distances}
    \label{sec-lb-distances}
    Having established all upper bounds in Theorem~\ref{thrm:fpp}, we continue with the lower bounds. We analyse the lower bound on the typical distance in \eqref{typical} in Section \ref{sec-typical}, on the flooding in \eqref{flooding} in Section \ref{sec:floodlb} and on the diameter in \eqref{diameter} in Section \ref{sec-diameter-lb}.
	
	\subsection{Lower bounds on typical distances}
	\label{sec-typical}
    Let us start by explaining how we proceed. To analyse the distance between two vertices, selected uniformly at random from $\T_n^d$, as in \eqref{typical}, we show that their exploration processes are disjoint up to time $(1-\eps)\frac12 R_n^{-1}\log n$ with high probability. This implies that their typical distance is at least $(1-\eps)R_n^{-1}\log n$ with high probability. 

    Let $U$ be vertex selected uniformly at random from $\T_n^d$ and let $v\in \T_n^d$ fixed. First, we write 
    \be \ba\label{eq:firstbound}
    \P{X_{U,v}\leq 2t_n}&\leq \P{U=v}+\P{B_U(t_n)\cap  B_v(t_n)\neq \emptyset\,|\, U\neq v}\\ 
    &=\frac1n +\P{\exists w\in \T_n^d: w\in B_U(t_n), w\in B_v(t_n)\,|\, U\neq v}\\
    &\leq \frac1n +\sum_{w\in \T_n^d}\P{X_{U,w}\leq t_n, \wt X_{v,w}\leq t_n\,|\, U\neq v},
    \ea\ee 
    where $\wt X_{v,w}$ is the time for an exploration process started from $v$ to reach $w$ outside of $B_U(t_n)$, i.e.\ on $\T_n^d\backslash B_U(t_n)$. Since $\wt X_{v,w}\succeq X_{v,w}'$, where $X_{v,w}$ is a copy of $X_{v,w}$, independent of $X_{U,w}$, we obtain the upper bound
    \be\label{eq:secondbound}
    \frac1n +\sum_{w\in \T_n^d}\P{X_{U,w}\leq t_n\,|\, U\neq v}\P{X_{v,w}\leq t_n}\leq \frac1n+ \sum_{w\in \T_n^d}\frac{1}{n-1}\E{|B_w(t_n)|}\P{X_{v,w}\leq t_n}.
    \ee 
    By translation invariance, $\E{|B_w(t_n)|}=\E{|B_{\mathbf 0}(t_n)|}$, where $\mathbf 0$ denotes the origin. We thus arrive at 
    \be 
    \frac1n+ \frac{1}{n-1}\E{|B_{\mathbf 0}(t_n)|}\sum_{w\in \T_n^d}\P{X_{v,w}\leq t_n}=\frac1n +\frac{1}{n-1}\E{|B_{\mathbf 0}(t_n)|}^2.
    \ee 
    Finally, we use that by the upper bound in~\eqref{eq:Lambdabounds}, 
    \be \label{eq:expballsize}
    \E{|B_{\mathbf 0}(t_n)|}^2\leq \big(\e^{R_nt_n}\big)^2=n^{1-\eps}, 
    \ee 
    to arrive at the desired result that 
    \be 
    \P{X_{U,v}\leq (1-\eps)R_n^{-1}\log n}=\P{X_{U,v}\leq 2t_n}=o(1).
    \ee 
    An almost analogous argument can be used to yield the same result for two uniform vertices $U$ and $V$ as well. This proves the lower bound for~\eqref{typical} in Theorem~\ref{thrm:fpp}.\qed

    \subsection{Lower bound on flooding time}
	\label{sec:floodlb}

    In this section we prove a lower bound on the flooding time, and is an adaptation of the idea for a lower bound on the \emph{diameter} for the case $\alpha=0$ by Janson in~\cite{Jans99}. 
    
    Fix $v\in \T_n^d$ and $\eps>0$ small, and take some set $A\subset \T_n^d\backslash\{v\}$ such that $n_A:=|A|=\lfloor n^{1-\eps}\rfloor$ and let $B:=\T_n^d\backslash A$. We enumerate the vertices in $A$ in an arbitrary order, so that we can write $A=\{u_1,\ldots, u_{n_A}\}$. Abusing notation, we write $i\in[n_A]$ when we mean $u_i\in A$. 
    
    We define, for $i\in [n_A]$, the random variable $U_i:=\min_{v\in B}T_{u_iv}$. That is, $U_i$ is the transmission time from $u_i$ to $B$. We also set $b=(1-2\eps)R_n^{-1}\log n$ and $\mu_i:=\sum_{v\in B}\|u_i-v\|^{-\alpha}, i\in [n_A]$. For all $i\in[n_A]$, it readily holds that $\mu_i\leq R_n$. Furthermore,  independently of $i\in[n_A]$, since $\alpha<d$,
	\be \label{eq:mui}
	\mu_i=R_n-\sum_{\substack{v\in A\\ v\neq u_i}}\|u_i-v\|^{-\alpha}\geq R_n-R_{n_A}=R_n(1-\cO(n^{\eps(1-\alpha/d)})), 
	\ee 
    so that $\mu_i=R_n(1-\cO(n^{\eps(1-\alpha/d)}))$, where the $\cO$ term is independent of $i\in[n_A]$.
    Now, for $i\in[n_A]$, define
	\be 
	f_i(x):=-\frac1\mu_i\log\big(\e^{-b\mu_i}+(1-\e^{-b\mu_i}\big)\e^{-\mu_i x}\big). 
	\ee 
	By~\cite[Lemma $2.1$]{Jans99}, it follows that $f_i(U_i)\overset d= (U_i|U_i\leq b)$, i.e.\ $f_i(U_i)$ has the same distribution as $U_i$, conditionally on the event $U_i\leq b$. 
	
	Now, fix some $k\in [n_A]$ and let
	\be 
	U_i':=\begin{cases}
		f_i(U_i)&\mbox{if }i<k,\\
		U_i+b&\mbox{if }i=k, \\
		U_i&\mbox{if }i>k.
	\end{cases}
	\ee 
	If we let $\mathcal E_k:=\{U_k> b, U_i\leq b$ for all $i<k\}$, then, conditionally on $\mathcal E_k$, $U_i'$ has the same distribution as $U_i$, by the memoryless property of the exponential distribution. As the $U_i$ are independent, the distributional equality also holds jointly. Similarly, we set 
	\be 
	T_{u_iv}':=\begin{cases}
		T_{u_iv} -U_i+U_i'&\mbox{if }i\in [n_A], v\in B,\\
		T_{u_iv}&\mbox{otherwise.}
	\end{cases}
	\ee 
	Equivalently, $T_{vu_i}'=T_{u_iv}'$.	Again, by the memoryless property and since $T_{u_iv}-U_i$ is independent of $U_i$, it follows that, conditionally on $\mathcal E_k$, $T_{u_iv}'$ has the same distribution as $T_{u_iv}$. In an equivalent manner we define $X'_{v,u}$ for $v,u\in \T_n^d$, which, for each pair $v,u\in \T_n^d$, we can think of as  equal in distribution to $X_{v,u}$, conditionally on $\cE_k$.
    
    Suppose the following holds:
    \be\ba \label{eq:Tboundrad}
	U_i'&\geq (1-2\eps)U_i\qquad&&\text{for all }i\in [n_A],\\
	T_{u_ki}&\geq 3R_n^{-1}\log n\qquad &&\text{for all }i\in [n_A]\backslash \{k\},\\
	X_{u_k,v}&\geq (1-2\eps)R_n^{-1}\log n.
	\ea \ee 
    First, we observe that 
	\be \label{eq:Tbound2}
	T_{uv}'\geq (1-2\eps)T_{uv}\qquad \text{for all }v\neq u.
	\ee 
	To see this, we need only consider $u\in A, v\in B$. Suppose $u=u_i$ for some $i\in [n_A]$. Then, by the first line of~\eqref{eq:Tboundrad},
	\be 
	T_{u_iv}'=T_{u_iv}-U_i+U_i'\geq T_{u_iv}-2\eps U_i\geq (1-2\eps)T_{u_iv},
	\ee 
	as $U_i$ is the minimal edge-weight between $u_i$ and $B$. As $X_{u_k,v}\geq(1-2\eps)R_n^{-1}\log n$, it implies that for any path $(u_k=u_{i_0},u_{i_1},\ldots, u_{i_\ell}=v)$ it holds that
    \be 
    W=\sum_{j=1}^\ell T_{u_{j-1}u_j}\geq (1-2\eps)R_n^{-1}\log n. 
    \ee 
    We then set $W':=\sum_{j=1}^\ell T_{u_{j-1}u_j}'$ to be the conditional weight of this path. If $u_{i_1}\in A$, then $W'\geq T_{u_ku_{i_1}}\geq 2R_n^{-1}\log n$ by the second line of~\eqref{eq:Tboundrad}. If instead $u_{i_1}\not\in A$, then $T_{u_ku_{i_1}}'=T_{u_ku_{i_1}}+b$ and 
    \be 
    W'\geq b+(1-2\eps)W\geq (1-2\eps)R_n^{-1}\log n+(1-2\eps)^2R_n^{-1}\log n \geq (2-6\eps)R_n^{-1}\log n.
    \ee 
    It thus follows that, for $k\in[n_A]$,
    \be 
    \P{X_{v,u_k}\geq(2-6\eps)R_n^{-1}\log n\,|\, \cE_k}=\P{X_{v,u_k}'\geq(2-6\eps)R_n^{-1}\log n}\geq \P{\eqref{eq:Tboundrad}\text{ holds for }k}.
    \ee 
    We can then write 
    \be \ba
    \P{\max_{u\in \T_n^d}X_{v,u}\geq (2-6\eps)R_n^{-1}\log n}&\geq \P{\bigcup_{u\in A}\{X_{v,u}\geq (2-6\eps)R_n^{-1}\log n\}}\\ 
    &\geq \P{\bigcup_{k\in [n_A]}\{X_{v,u_k}\geq(2-6\eps)R_n^{-1}\log n\}\cap \cE_k}.
    \ea \ee 
    Since the events $\cE_k$ are disjoint, we can write the union as a sum and introduce $N$ to be a uniform element of $[n_A]$, so that $u_N$ is a uniform element of $A$. Then, we obtain 
    \be \ba\label{eq:maxlb}
    n_A\P{\{X_{v,u_N}\geq (2-6\eps)R_n^{-1}\log n\}\cap \cE_N}&=n_A \P{X_{v,u_N}\geq (2-6\eps)R_n^{-1}\log n\,|\, \cE_N}\P{\cE_N}\\ 
    &=\P{X_{v,u_N}\geq (2-6\eps)R_n^{-1}\log n\,|\, \cE_N}\sum_{k=1}^{n_A}\P{\cE_k}\\
    &\geq \P{\eqref{eq:Tboundrad}\text{ holds for }N} \sum_{k=1}^{n_A}\P{\cE_k}.
    \ea \ee 
    Let us first argue that the sum is $1-o(1)$. As the events $\EE_k$ are disjoint and the random variables $(U_i)_{i\in A}$ independent, 
	\be 
	\sum_{k=1}^{n_A}\P{\EE_k}=\P{\bigcup_{k=1}^{n_A}\EE_k}=1-\P{\bigcap_{k=1}^{n_A}\EE_k^c}\geq 1-\P{U_i\leq b\text{ for all } i\in [n_A]}
	= 1-\prod_{i=1}^{n_A} \P{U_i\leq b}.
	\ee 
	Now, as $b=(1-2\eps)R_n^{-1}\log n$, $U_i\sim \text{Exp}(\mu_i)$, and $\mu_i=R_n(1-\cO(n^{-\eps(1-\alpha/d)}))$ by~\eqref{eq:mui}, where the $\cO(n^{-\eps(1-\alpha/d)})$ term is independent of $i$, we obtain 
	\be \label{eq:implied2}
	1-\prod_{i=1}^{n_A} \P{U_i\leq b}=1-\prod_{i=1}^{n_A}\big(1-\exp(-\mu_i b)\big)\geq 1-\exp(-n^{(2\eps-1)+o(1)}n_A)=1-\exp(-n^{\eps+o(1)}).
	\ee 
    We use this in~\eqref{eq:maxlb} to arrive at
    \be 
    \P{\max_{u\in \T_n^d }X_{v,u}\geq (2-6\eps)R_n^{-1}\log n}\geq (1-o(1))\P{\eqref{eq:Tboundrad}\text{ holds for }N}.
    \ee 
    To show that the probability on the right-hand side converges to one, it suffices to show that the events in~\eqref{eq:Tboundrad} hold with high probability, uniformly in $k\in A$, which we do now.

    Fix $k\in [n_A]$. For the {\bf first line} in~\eqref{eq:Tboundrad}, we have 
	\be 
	\P{\bigcup_{i\in A}\{U_i'<(1-2\eps)U_i\}}\leq\sum_{i\in A}\P{U_i'<(1-2\eps)U_i}=\sum_{i<k}\P{f_i(U_i)\leq (1-2\eps)U_i},
	\ee  
    where the last equality follows since $U_i'\geq U_i$ for all $i\geq k$ by \eqref{eq:Tboundrad}.
    
	It now follows from~\cite[Lemma $2.1$]{Jans99} and~\eqref{eq:mui} that  the right-hand side can be upper bounded by 
	\be 
	\sum_{i<k}\P{U_i>b/(1-2\eps)-1/\mu_i}=\sum_{i<k}\e^{1-b\mu_i/(1-2\eps)}\leq n_A\e^{1-bR_n(1-\cO(n^{-\eps}))/(1-2\eps)}=o(1), 
	\ee 
	since $bR_n/(1-2\eps)=\log n$ and $n_A=o(n)$, and where the right-hand side is independent of the choice of $k\in A$. As a result, the first line in~\eqref{eq:Tboundrad} also holds for a uniform element $N\in[n_A]$.
 
    For the {\bf second line} in~\eqref{eq:Tboundrad}, with $k\in A$ fixed,
	\be 
	\P{\bigcup_{i\in A\backslash\{k\}}\!\!\{T_{u_ku_i}<3R_n^{-1}\log n\}}=1-\exp\Big(-3\frac{\log n}{R_n}\sum_{i\in A\backslash\{k\}}\|u_k-u_i\|^{-\alpha}\Big)\leq 3\frac{\log n}{R_n}R_{n_A-1}=o(1),
	\ee 
    where the penultimate step uses the same idea as in~\eqref{bounds-Gamma-j} and where the upper bound holds uniformly in $k$ and thus also for $k=N$. 
    
    Finally, the {\bf third line} of~\eqref{eq:Tboundrad} requires that the typical distance between a uniform element $u_N$ from $A$ and $v$ is at least $(1-2\eps)R_n^{-1}\log n$ with high probability. The lower bound on the flooding time in Section~\ref{sec:floodlb} tells us this holds if we had a uniform vertex from $\T_n^d$ rather than from $A$. However, this argument is readily adapted to work for $u_N$ as well. Namely, using $u_N$ instead of $U$ in~\eqref{eq:firstbound} and~\eqref{eq:secondbound} (where we don't have to condition on $u_N\neq v$ since $v\not\in A$) yields for $t_n:=(1-2\eps)\frac12 R_n^{-1}\log n$,
    \be 
    \P{X_{u_N,v}\leq 2t_n}\leq \sum_{w\in \T_n^d}\P{X_{u_N,w}\leq t_n}\P{X_{v,w}\leq t_n}=\sum_{w\in \T_n^d}\P{X_{v,w}\leq t_n}\frac{1}{n_A} \E{|B_w(t_n)\cap A|}.
    \ee 
    We can simply bound $\E{|B_w(t_n)\cap A|}\leq \E{|B_w(t_n)|}=\E{|B_{\mathbf 0}(t_n)|}$, where $\mathbf 0$ is the origin. As a result, we obtain the upper bound 
    \be 
    \frac{1}{n_A} \E{|B_{\mathbf 0}(t_n)|}\sum_{w\in \T_n^d}\P{X_{v,w}\leq t_n}=\frac{1}{n_A} \E{|B_{\mathbf 0}(t_n)|}^2\leq \frac{1}{n_A}n^{1-2\eps}=o(1),
    \ee 
    where we use~\eqref{eq:expballsize} in the penultimate step. This proves the third line of~\eqref{eq:Tboundrad} and concludes the proof of the lower bound for the flooding time.
    \qed
    
	\subsection{Lower bound on the  diameter}
	\label{sec-diameter-lb}
	
	This section proves a lower bound on the diameter, and uses a similar approach and the same definitions as in Section~\ref{sec:floodlb}.
    
    Fix $k\in[n_A]$ and suppose that the first two lines of \eqref{eq:Tboundrad} hold, while the third line is replaced with
    \be\ba \label{eq:Tbound}
	Y_k:=\max_{\substack{v\in \T_n^d\\ v\neq u_k}}X_{u_k,v}&\geq (2-\eps)R_n^{-1}\log n.
    \ea \ee 
   Indeed, the statement in \eqref{eq:Tbound} holds with high probability by what is discussed in Section~\ref{sec:floodlb}, while \eqref{eq:Tboundrad} was already proved to hold with high probability.
    
    We then claim that it follows from the first two lines of~\eqref{eq:Tboundrad} and~\eqref{eq:Tbound} that $$Y_k':=\max_{\substack{v\in \T_n^d\\ v\neq u_k}}X_{u_k,v}'\geq (3-7\eps)R_n^{-1}\log n.$$ 
    Since $Y_k\geq (2-\eps)R_n^{-1}\log n$ by~\eqref{eq:Tbound}, there exists a path $(u_k=u_{i_0},u_{i_1},\ldots, u_{i_\ell})$, such that 
	\be \label{eq:Wbound}
	W:=\sum_{j=1}^\ell T_{u_{i_{j-1}}u_{i_j}}\geq (2-\eps)R_n^{-1}\log n. 
	\ee 
	Consider such a path, and its conditional weight $W'=\sum_{j=1}^\ell T'_{u_{i_{j-1}}u_{i_j}}$. If $i_1\in [n_A]$, then $W'\geq T_{u_ku_{i_1}}\geq 3R_n^{-1}\log n$ by the second line in~\eqref{eq:Tboundrad}. If $i_1\not\in [n_A]$, then $T'_{u_ku_{i_1}}=T_{u_ku_{i_1}}+b$, and 
	\be 
	W'\geq b+(1-2\eps)W\geq (1-2\eps)R_n^{-1}\log n+(1-2\eps)(2-\eps)R_n^{-1}\log n\geq (3-7\eps)R_n^{-1}\log n, 
	\ee
	where the second and third step follow from~\eqref{eq:Tbound2} and~\eqref{eq:Wbound}, respectively. It thus follows that this lower bound on the conditional transmission time holds for \emph{any} path from $u_k$ to $u_{i_\ell}$, so that $Y'_k\geq X'_{u_ku_{i_\ell}}\geq (3-7\eps)R_n^{-1}\log n$. If we let $\cC_{k,n}$ denote the event that the first two lines of~\eqref{eq:Tboundrad} and~\eqref{eq:Tbound} hold, we have thus shown that
	\be \label{eq:implied}
	\P{Y_k\geq (3-7\eps)R_n^{-1}\log n\,|\, \mathcal E_k}=\P{Y_k'\geq (3-7\eps)R_n^{-1}\log n}\geq \P{\cC_{k,n}}.
	\ee 
	We next use the fact that the events $(\cE_k)_{k\in[n_A]}$ are disjoint to obtain the lower bound
	\be \ba\label{eq:diamlb}
	\P{\max_{u,v\in \T_n^d}X_{u,v}\geq (3-7\eps)R_n^{-1}\log n}&\geq \sum_{k=1}^{n_A}\P{Y_k\geq (3-7\eps)R_n^{-1}\log n\,|\, \mathcal E_k}\P{\mathcal E_k}\\
	&\geq \min_{k\in [n_A]}\P{Y_k\geq (3-7\eps)R_n^{-1}\log n \,|\, \mathcal E_k}\sum_{k=1}^{n_A}\P{\EE_k}\\
    &=(1-o(1))\min_{k\in [n_A]}\P{\cC_{k,n}},
	\ea \ee 
    where the final step follows from~\eqref{eq:implied2}. Since we have already verified  that the event $\cC_{k,n}$ holds with high probability independently of the choice of $k\in A$, we conclude that the right-hand side converges to 1 as $n\rightarrow \infty$. This proves the lower bound for~\eqref{diameter} in Theorem~\ref{thrm:fpp}.
	\qed

	\appendix 
	\section{The limiting constants}\label{sec:app}
	In this appendix we investigate the limiting constant
	\be 
	R=R(d,p,\alpha):=\lim_{n\to\infty} R_n/n^{1-\alpha/d},  
	\ee 
    when $\alpha<d$, where we recall $R_n$ from~\eqref{eq:rn}. By switching from summation to integration and using a variable transform $x=n^{1/d}w$, we obtain
	\be \label{eq:rn-b}
	R_n=(1+o(1))\int_{[-\frac12 n^{1/d},\frac12 n^{1/d}]^d}\|x\|_p^{-\alpha}\,\d x=(1+o(1))n^{1-\alpha/d}\int_{[-\frac12,\frac12]^d}\|w\|_p^{-\alpha}\,\d w.
	\ee
	Using symmetry and again using a variable transform $y=2w$, we arrive at
	\be 
	R_n=(1+o(1))2^dn^{1-\alpha/d}\int_{[0,\frac12]^d}\|w\|_p^{-\alpha}\,\d w=(1+o(1))2^\alpha n^{1-\alpha/d}\int_{[0,1]^d}\|y\|_p^{-\alpha}\,\d y.
	\ee 
	It follows that 
	\be\label{eq:Rint}
	R(d,p,\alpha)=2^\alpha\int_{[0,1]^d} \|y\|_p^{-\alpha}\,\d y.
	\ee
	Then, it is immediate that $R(d,p,0)=1$ for all $d\in\N$ and $p\in[1,\infty)\cup\{\infty\}$. It thus remains to consider $\alpha\in (0,d)$. For $p=\infty$, we can use that 
	\be \label{eq:gammafunc}
	\int_0^{\infty} t^{\alpha-1} \e^{-at}\,\d t=a^{-\alpha} \Gamma(\alpha),
	\ee 
	to write
	\be\ba \label{eq:Iinfty}
	R(d,\infty,\alpha)&=\frac{2^\alpha}{\Gamma(\alpha)}\int_{[0,1]^d} \int_0^{\infty} t^{\alpha-1} \e^{-t\max\{y_1,\ldots,y_d\}}\,\d t \d y_1\cdots \d y_d\\
	&=\frac{2^\alpha}{\Gamma(\alpha)}\int_0^{\infty} t^{\alpha-1} \int_{[0,1]^d} \e^{-t\max\{y_1,\ldots, y_d\}}\,\d y_1\cdots \d y_d \d t\\
	&=\frac{2^\alpha d!}{\Gamma(\alpha)}\int_0^{\infty} t^{\alpha-1} \int_0^1 \int_{y_d}^1\cdots \int_{y_2}^1 \e^{-ty_1}\,\d y_1\cdots \d y_d \d t.
	\ea \ee 
	Let us define the inner integrals as 
	\be 
	I(d,t):= \int_0^1 \int_{y_d}^1\cdots \int_{y_2}^1 \e^{-ty_1}\,\d y_1\cdots \d y_d.
	\ee 
	By evaluating the innermost integral and observing that $I(1,t)=1-\e^{-t}$, we obtain
	\be 
	I(d,t)=\frac1t I(d-1,t)-\frac1t\e^{-t}\frac{1}{(d-1)!}=\ldots =t^{-d}-\sum_{i=1}^d t^{-i}\e^{-t}\frac{1}{(d-i)!}=t^{-d}\bigg(1-\sum_{j=0}^{d-1}\e^{-t}\frac{t^j}{j!}\bigg).
	\ee 
	The part between brackets equals the probability that a rate $t$ Poisson random variable is at least $d$. Via a rate-one Poisson process we can see that this probability equals the probability that a Gamma$(d,1)$ random variable is at most $t$. Using this in~\eqref{eq:Iinfty}, we arrive at 
	\be \ba\label{eq:Rinfty}
	R(d,\infty,\alpha)&=\frac{2^\alpha d!}{\Gamma(\alpha)}\int_0^\infty t^{\alpha-d-1}\int_0^t \frac{1}{(d-1)!}x^{d-1}\e^{-x}\,\d x\d t\\
	&=\frac{2^\alpha d}{\Gamma(\alpha)}\int_0^\infty x^{d-1}\e^{-x}\int_x^\infty t^{\alpha-d-1}\,\d t \d x\\
	&=\frac{2^\alpha}{\Gamma(\alpha)}\frac{d}{d-\alpha}\int_0^\infty x^{\alpha-1}\e^{-x}\,\d x=\frac{d}{d-\alpha} 2^\alpha.
	\ea \ee 
	For $p\in[1,\infty)$ a similar approach yields
	\eqan{
		R(d,p,\alpha)&=2^\alpha \int_{[0,1]^d} \Big(\sum_{i=1}^d y_i^p\Big)^{-\alpha/p}\,\d y_1\cdots \d y_d\\
		&=\frac{2^\alpha}{\Gamma(\alpha/p)}\int_{[0,1]^d}\int_0^\infty t^{\alpha/p-1}\exp\Big(-t\sum_{i=1}^d y_i^p\Big)\,\d t\d y_1\cdots \d y_d\\
		&=\frac{2^\alpha}{\Gamma(\alpha/p)}\int_0^\infty t^{\alpha/p-1}\int_{[0,1]^d}\exp\Big(-t\sum_{i=1}^d y_i^p\Big)\,\d y_1\cdots \d y_d\d t.
	}
	We use the variable transformation $z_i=ty_i^p$ to then obtain
	\eqan{
    \label{eq:finint}
	R(d,p,\alpha)&=\frac{2^\alpha}{\Gamma(\alpha/p)}\int_0^\infty t^{\alpha/p-1}\bigg(\int_0^t \frac1p t^{-1/p}z^{1/p-1}\e^{-z}\,\d z\bigg)^d\d t\\
    &=\frac{2^\alpha p^{-d}}{\Gamma(\alpha/p)}\int_0^\infty t^{(\alpha-d)/p-1}\gamma(1/p,t)^d\,\d t, \nn
	}
	where $\gamma(1/p,t):=\int_0^t z^{1/p-1}\e^{-z}\,\d z$ is the lower incomplete gamma function. 
	
	We are unable to explicitly evaluate this integral in full generality. For $d=2$ and $\alpha\in(0,2)$, however, an explicit expression of $R$ in terms of the hypergeometric function can be obtained for all $p\in[1,\infty)$. Namely, the integral in~\eqref{eq:finint} can be written, using integration by parts, as
	\be \ba\label{eq:intpart}
	\Big[{}&\frac{p}{\alpha-2}t^{(\alpha-2)/p}\gamma(1/p,t)^2\Big]_0^\infty +\frac{2p}{2-\alpha}\int_0^\infty t^{(\alpha-2)/p}\gamma(1/p,t)t^{1/p-1}\e^{-t}\,\d t\\
	&=\frac{2p}{2-\alpha}\int_0^\infty t^{(\alpha-1)/p-1}\gamma(1/p,t)\e^{-t}\,\d t,
	\ea \ee 
	where the last step uses that $\lim_{t\downarrow 0}t^{-1/p}\gamma(1/p,t)=p$ and $\alpha\in(0,2)$. Now, Gradshteyn and Ryzhik~\cite[6.455.2]{GraRyz65} show this equals 
	\be 
	\frac{2p^2}{2-\alpha}\Gamma(\alpha/p)2^{-\alpha/p} {}_2F_1\Big(1,\frac{\alpha}{p},1+\frac1p,\frac12\Big),
	\ee 
	where $_2F_1(a,b,c,z)$ denotes the hypergeometric function, defined for $|z|<1$ as
	\be 
	_2F_1(a,b,c,z):=\sum_{n=0}^\infty \frac{\Gamma(a+n)\Gamma(c+n)\Gamma(b)}{\Gamma(a)\Gamma(c)\Gamma(b+n)}\frac{z^n}{n!}.
	\ee 
	Substituting this into~\eqref{eq:finint} finally yields 
	\be \label{eq:Rp}
	R(2,p,\alpha)=\frac{2^{1+\alpha(1-1/p)}}{2-\alpha}{}_2F_1\Big(1,\frac{\alpha}{p},1+\frac1p,\frac12\Big).
	\ee 
    We also observe that the limit of $p\to \infty$ on the right-hand side yields the same expression as for $R(d,\infty,\alpha)$ in~\eqref{eq:Rinfty}. Since the supremum-norm corresponds to the limit of the $p$-norm as $p\to\infty$, it follows that~\eqref{eq:Rinfty} also follows from (the derivation of)~\eqref{eq:Rp} and the dominated convergence theorem.
    
    In the case of general dimension $d$, one interesting observation is the following: let $M_{d,p}$ be the maximum of $d$ i.i.d.\ Gamma$(1/p,1)$ random variables and let $f_{M_d}$ denote its probability density function. Then, 
	\be 
	\gamma(1/p,t)^d=\Gamma(1/p)^d\P{M_{d,p}\leq t}.
	\ee 
	As a result, using integration by parts on the right-hand side of~\eqref{eq:finint} (as in~\eqref{eq:intpart}) we obtain
	\be \label{eq:Rexp}
	R(d,p,\alpha)=\frac{2^\alpha\Gamma(1/p)^d p^{-(d-1)}}{\Gamma(\alpha/p)(d-\alpha)}\int_0^\infty t^{(\alpha-d)/p}f_{M_{d,p}}(t)\,\d t=\frac{2^\alpha\Gamma(1/p)^d p^{-(d-1)}}{\Gamma(\alpha/p)(d-\alpha)}\E{M_{d,p}^{(\alpha-d)/p}}.
	\ee 
 Using that the density $f_{M_{d,p}}$ is given by
    \[
    f_{M_{d,p}}(t)=dF_{\Gamma}(t)^{d-1}f_{\Gamma}(t),
    \]
where $f_{\Gamma}$ is the density of a Gamma$(1/p,1)$ random variable, we can rewrite
    \eqan{
    \E{M_{d,p}^{(\alpha-d)/p}}
    &=d \int_0^\infty t^{(\alpha-d)/p}
    F_{\Gamma}(t)^{d-1}f_{\Gamma}(t)\,\d t\nn\\
    &=\frac{d\Gamma((\alpha-(d-1))/p)}
    {\Gamma(1/p)}\prob(\Gamma'>\max_{i=1}^{d-1}\Gamma_i),\nn
    }
where $\Gamma'$ is a Gamma$((\alpha-(d-1))/p,1)$ random variable, and $(\Gamma_i)_{i=1}^{d-1}$ are i.i.d.\ Gamma$(1/p,1)$ random variables. While these provide analytic and probabilistic expressions for $R(d,p,\alpha)$, they do not allow us to simplify the formulas further.

	\medskip

	\paragraph{\bf Acknowledgements.} This work initiated during the RandNET Workshop on Random Graphs held 22-30 August 2022 in Eindhoven, the Netherlands. BL would like to thank Marcos Kiwi, Lyuben Lichev, Amitai Linker, and Dieter Mitsche for discussions on a closely related problem that served as inspiration for the open problem he presented at the aforementioned workshop. 
 
    The work of RvdH was supported in part by the Netherlands Organisation for Scientific Research (NWO) through Gravitation-grant {\sc NETWORKS}-024.002.003. BL was supported through the grant GrHyDy ANR-20-CE40-0002.

    The authors would like to thank an anonymous reviewer for a careful reading of the manuscript, and for their suggestions regarding the introduction and literature review.
	
	\bibliographystyle{abbrv}
	\bibliography{distfpp}

\providecommand{\vander}{van der }
\begin{thebibliography}{10}

\bibitem{AddBrouLug10}
L.~Addario-Berry, N.~Broutin, and G.~Lugosi.
\newblock The longest minimum-weight path in a complete graph.
\newblock {\em Combinatorics, Probability and Computing}, 19(1):1–19, 2010.

\bibitem{AdrKom18}
E.~Adriaans and J.~Komj{\'a}thy.
\newblock Weighted distances in scale-free configuration models.
\newblock {\em Journal of Statistical Physics}, 173:1082--1109, 2018.

\bibitem{BarHofKom17}
E.~Baroni, R.~van~der Hofstad, and J.~Komj{\'a}thy.
\newblock Nonuniversality of weighted random graphs with infinite variance
  degree.
\newblock {\em Journal of Applied Probability}, 54(1):146--164, 2017.

\bibitem{BarHofKom19}
E.~Baroni, R.~van~der Hofstad, and J.~Komj{\'a}thy.
\newblock Tight fluctuations of weight-distances in random graphs with
  infinite-variance degrees.
\newblock {\em Journal of Statistical Physics}, 174:906--934, 2019.

\bibitem{BenKesPerSch11}
I.~Benjamini, H.~Kesten, Y.~Peres, and O.~Schramm.
\newblock Geometry of the uniform spanning forest: transitions in dimensions 4,
  8, 12,….
\newblock {\em Selected Works of Oded Schramm}, pages 751--777, 2011.

\bibitem{BhaHof12}
S.~Bhamidi and R.~{\vander H}ofstad.
\newblock Weak disorder asymptotics in the stochastic mean-field model of
  distance.
\newblock {\em The Annals of Applied Probability}, 22(1):29--69, 2012.

\bibitem{BhaHof17}
S.~Bhamidi and R.~{\vander H}ofstad.
\newblock Diameter of the stochastic mean-field model of distance.
\newblock {\em Combinatorics, Probability and Computing}, 26(6):797--825, 2017.

\bibitem{BhaHofHoog13}
S.~Bhamidi, R.~{\vander H}ofstad, and G.~Hooghiemstra.
\newblock Weak disorder in the stochastic mean-field model of distance {II}.
\newblock {\em Bernoulli}, 19(2):363--386, 2013.

\bibitem{BhaHofHoog10.2}
S.~Bhamidi, R.~van~der Hofstad, and G.~Hooghiemstra.
\newblock {First passage percolation on random graphs with finite mean
  degrees}.
\newblock {\em The Annals of Applied Probability}, 20(5):1907 -- 1965, 2010.

\bibitem{BhaHofHoog11}
S.~Bhamidi, R.~van~der Hofstad, and G.~Hooghiemstra.
\newblock First passage percolation on the {E}rd{\H{o}}s--{R}{\'e}nyi random
  graph.
\newblock {\em Combinatorics, Probability and Computing}, 20(5):683--707, 2011.

\bibitem{BhaHofHoog17}
S.~Bhamidi, R.~van~der Hofstad, and G.~Hooghiemstra.
\newblock {Universality for first passage percolation on sparse random graphs}.
\newblock {\em The Annals of Probability}, 45(4):2568 -- 2630, 2017.

\bibitem{BisKrie24}
M.~Biskup and A.~Krieger.
\newblock Arithmetic oscillations of the chemical distance in long-range
  percolation on z d.
\newblock {\em The Annals of Applied Probability}, 34(3):2986--3017, 2024.

\bibitem{BroaHam57}
S.~R. Broadbent and J.~M. Hammersley.
\newblock Percolation processes: I. crystals and mazes.
\newblock In {\em Mathematical proceedings of the Cambridge philosophical
  society}, volume~53, pages 629--641. Cambridge University Press, 1957.

\bibitem{ChaDey16}
S.~Chatterjee and P.~S.~Dey.
\newblock Multiple phase transitions in long-range first-passage percolation on
  square lattices.
\newblock {\em Communications on Pure and Applied Mathematics}, 69(2):203--256,
  2016.

\bibitem{EckGoodHofNar13}
M.~Eckhoff, J.~Goodman, R.~{\vander H}ofstad, and F.~R. Nardi.
\newblock Short paths for first passage percolation on the complete graph.
\newblock {\em Journal of Statistical Physics}, 151(6):1056--1088, 2013.

\bibitem{EckGoodHofNar15.1}
M.~Eckhoff, J.~Goodman, R.~van~der Hofstad, and F.~R. Nardi.
\newblock {Long paths in first passage percolation on the complete graph {I}.
  Local {PWIT} dynamics}.
\newblock {\em Electronic Journal of Probability}, 25:1 -- 45, 2020.

\bibitem{EckGoodHofNar15.2}
M.~Eckhoff, J.~Goodman, R.~van~der Hofstad, and F.~R. Nardi.
\newblock Long paths in first passage percolation on the complete graph {II}.
  {G}lobal branching dynamics.
\newblock {\em Journal of Statistical Physics}, 181(2):364--447, 2020.

\bibitem{GraRyz65}
I.~S. Gradshteyn and I.~M. Ryzhik.
\newblock {\em Table of integrals, series, and products}.
\newblock Fourth edition prepared by Ju. V. Geronimus and M. Ju. Ce\u\i tlin.
  Translated from the Russian by Scripta Technica, Inc. Translation edited by
  Alan Jeffrey. Academic Press, New York, 1965.

\bibitem{HofHoogMieg02}
R.~{\vander H}ofstad, G.~Hooghiemstra, and P.~Van~Mieghem.
\newblock The flooding time in random graphs.
\newblock {\em Extremes}, 5(2):111--129, 2002.

\bibitem{HofKom17}
R.~{\vander H}ofstad and J.~Komj{\'a}thy.
\newblock Explosion and distances in scale-free percolation.
\newblock {\em arXiv preprint arXiv:1706.02597}, 2017.

\bibitem{HofHoogMie01}
R.~v.~d. Hofstad, G.~Hooghiemstra, and P.~Van~Mieghem.
\newblock First-passage percolation on the random graph.
\newblock {\em Probability in the Engineering and Informational Sciences},
  15(2):225–237, 2001.

\bibitem{HofHoogMie06}
R.~v.~d. Hofstad, G.~Hooghiemstra, and P.~Van~Mieghem.
\newblock Size and weight of shortest path trees with exponential link weights.
\newblock {\em Combinatorics, Probability and Computing}, 15(6):903–926,
  2006.

\bibitem{HoogMie08}
G.~Hooghiemstra and P.~Van~Mieghem.
\newblock The weight and hopcount of the shortest path in the complete graph
  with exponential weights.
\newblock {\em Combinatorics, Probability and Computing}, 17(4):537–548,
  2008.

\bibitem{Jans99}
S.~Janson.
\newblock One, two and three times log n/n for paths in a complete graph with
  random weights.
\newblock {\em Combinatorics, Probability and Computing}, 8(4):347--361, 1999.

\bibitem{JorKom20}
J.~Jorritsma and J.~Komj{\'a}thy.
\newblock Weighted distances in scale-free preferential attachment models.
\newblock {\em Random Structures \& Algorithms}, 57(3):823--859, 2020.

\bibitem{JorKom22}
J.~Jorritsma and J.~Komj{\'a}thy.
\newblock Distance evolutions in growing preferential attachment graphs.
\newblock {\em The Annals of Applied Probability}, 32(6):4356--4397, 2022.

\bibitem{KomLapLen21}
J.~Komj{\'a}thy, J.~Lapinskas, and J.~Lengler.
\newblock Penalising transmission to hubs in scale-free spatial random graphs.
\newblock In {\em Annales de l'Institut Henri Poincar{\'e} (B) Probabilit{\'e}s
  et statistiques}, volume~57, pages 1968--2016. Institut Henri Poincar{\'e},
  2021.

\bibitem{KomLapLenSchal23_1}
J.~Komj{\'a}thy, J.~Lapinskas, J.~Lengler, and U.~Schaller.
\newblock Four universal growth regimes in degree-dependent first passage
  percolation on spatial random graphs {I}.
\newblock {\em arXiv preprint arXiv:2309.11840}, 2023.

\bibitem{KomLapLenSchal23_2}
J.~Komj{\'a}thy, J.~Lapinskas, J.~Lengler, and U.~Schaller.
\newblock Four universal growth regimes in degree-dependent first passage
  percolation on spatial random graphs {II}.
\newblock {\em arXiv preprint arXiv:2309.11880}, 2023.

\bibitem{KomLod20}
J.~Komj{\'a}thy and B.~Lodewijks.
\newblock Explosion in weighted hyperbolic random graphs and geometric
  inhomogeneous random graphs.
\newblock {\em Stochastic Processes and their Applications}, 130(3):1309--1367,
  2020.

\end{thebibliography}

\end{document}